\documentclass[12pt]{article}
\usepackage{mathtools} 
\usepackage{amsmath, amssymb} 
\usepackage{amsthm}
\usepackage{enumerate}  
\usepackage{url} 
\usepackage{authblk}
\usepackage{tikz}
\usetikzlibrary{decorations.markings}
\usepackage{hyperref}

\usepackage{graphicx}
\usepackage{caption,subcaption}

\usepackage{multirow}

\newcommand\cx{{\mathbb C}}

\newcommand\re{{\mathbb R}}

\DeclarePairedDelimiter\abs{\lvert}{\rvert}%
\DeclarePairedDelimiter\norm{\lVert}{\rVert}%

\makeatletter
\let\oldabs\abs
\def\abs{\@ifstar{\oldabs}{\oldabs*}}
\let\oldnorm\norm
\def\norm{\@ifstar{\oldnorm}{\oldnorm*}}
\makeatother

\newcommand\comp[1]{{\mkern2mu\overline{\mkern-2mu#1}}}
\newcommand\pmat[1]{\begin{pmatrix} #1 \end{pmatrix}}
\newcommand\seq[4]{#1_{#2},#1_{#3},\ldots,#1_{#4}}

\newtheoremstyle{plainsl}%
	{\topsep}
	{\topsep}
	{\slshape} 
	{}
	{\normalfont\bfseries}
	{.}
	{ }
	{}

\swapnumbers

{\theoremstyle{plainsl}
\newtheorem{theorem}{Theorem}[section]
\newtheorem{lemma}[theorem]{Lemma}
\newtheorem{corollary}[theorem]{Corollary}}
{\theoremstyle{remark}
}

\renewcommand\proof{\noindent\textsl{Proof. }}
\newcommand\sqr[2]{{\vbox{\hrule height.#2pt
    \hbox{\vrule width.#2pt height#1pt \kern#1pt
        \vrule width.#2pt}\hrule height.#2pt}}}
\renewcommand\qed{%
	\ifmmode\eqno\sqr53
	\else\nolinebreak\ \hfill\sqr53\medbreak\fi}

\DeclareMathOperator{\col}{col}

\newcommand\ip[2]{\langle#1,#2\rangle}
\newcommand\one{{\bf1}}

\usepackage{blkarray}

\title{The average search probabilities of discrete-time quantum walks}

\author{Hanmeng Zhan}

\affil{Department of Mathematics and Statistics, York University, Toronto, ON, Canada\\\texttt{h3zhan@yorku.ca}}

\begin{document}
\maketitle

\begin{abstract}
We study the average probability that a discrete-time quantum walk finds a marked vertex on a graph. We first show that, for a regular graph, the spectrum of the transition matrix is determined by the weighted adjacency matrix of an augmented graph. We then consider the average search probability on a distance regular graph, and find a formula in terms of the adjacency matrix of its vertex-deleted subgraph. In particular, for any family of 
\begin{itemize}
\item complete graphs, or
\item strongly regular graphs, or
\item distance regular graphs of a fixed parameter $d$, varying valency $k$ and varying size $n$, such that $k^{d-1}/n$ vanishes as $k$ increases,
\end{itemize}
the average search probability approaches $1/4$ as the valency goes to infinity. We also present a more relaxed criterion, in terms of the intersection array, for this limit to be approached by distance regular graphs.
\end{abstract}

\section{Introduction}

Quantum walks can be turned into algorithms to search marked vertices in graphs. While this idea was formalized by Shenvi, Kempe and Whaley \cite{Shenvi2003}, its first application dates back to 1996, when Grover \cite{Grover1996} showed that finding a marked vertex in a looped $K_n$ takes $O(\sqrt{n})$ steps of a quantum walk. Since then, quantum walk search has been studied on various graphs, including hypercubes \cite{Shenvi2003}, Cartesian powers of cycles \cite{Childs2003}, strongly regular graphs \cite{Janmark2014}, certain Johnson graphs \cite{Ambainis2003,Wong2016, Xue2019, Tanaka2021}, and more generally, regular locally arc-transitive graphs \cite{Hoyer2020}. For most of these graphs, quantum walks arrive at the marked vertices faster than the classical random walks.

In this paper, we consider a related problem: given a graph and a marked vertex, what is the average probability, over any period of time, that a discrete-time quantum walk finds the marked vertex? This probability converges as the time period tends to infinity, and we will refer to its limit as the \textsl{average search probability}. Intuitively, the average search probability should depend on the parameters of the graph. However, as we will see later, for many parametric families of distance regular graphs with a fixed diameter, the average probability approaches $1/4$ as the valency increases, regardless of the defining parameters. This surprising phenomenon echoes with the fast search of quantum walks on highly regular graphs.

We will start with a general observation, on the spectrum of a quantum walk that incorporates an oracle (Theorem \ref{thm:spd}). We then compute the average search probability for a distance regular graph, which, by results from equitable partitions, is a function in the eigenprojections of the vertex-deleted subgraph  (Theorem \ref{thm:avgsearch}). Finally, we derive conditions for the average search probability to approach $1/4$ as the valency increases (Theorem \ref{thm:criterion}), and show that this happens for any family of 
\begin{itemize}
\item complete graphs, or
\item strongly regular graphs, or
\item distance regular graphs of a fixed parameter $d$, varying valency $k$ and varying size $n$, such that $k^{d-1}/n$ vanishes as $k$ tends to infinity.
\end{itemize}

\section{Quantum walks with oracles \label{sec:model}}

Let $X$ be a $k$-regular graph on $n$ vertices, and let $a$ be a marked vertex of $X$. Our goal is to find $a$ using a discrete-time quantum walk on $X$.

Throughout, we will view each edge $\{u,v\}$ of $X$ as a pair of arcs $(u,v)$ and $(v,u)$. The \textsl{states} associated with $X$ are complex-valued functions on the arcs; hence, they form a vector space isomorphic to $\cx^n \otimes \cx^k$. The initial state of our quantum walks is the unit vector
\[x_0 = \frac{1}{\sqrt{nk}} \one_n\otimes \one_k.\]
To find the marked vertex, we apply a unitary matrix $U$, called the \textsl{transition matrix}, iteratively to the initial state. More specifically, $U$ is the product of the following unitary operators on $\cx^n \otimes \cx^k$: the \textsl{arc-reversal matrix} $R$, which represents the permutation that swaps $(u,v)$ with $(v,u)$; the \textsl{coin matrix}
\[C = I_n \otimes \left(\frac{2}{k} J_k-I_k\right),\]
and the \textsl{oracle}
\[O_a= (2 E_{aa}-I_n)\otimes I_k.\]
For algorithmic meanings of these operators, see Shenvi, Kempe and Whaley \cite{Shenvi2003}.

Now set
\[U = RCO_a.\]
At time $t$, our quantum walk will be in state $U^t x_0$ if it started with state $x_0$. As $U^t x_0$ is a unit vector, the entrywise product
\[(U^t x_0) \circ \comp{(U^t x_0)}\]
represents a probability distribution, and its $(a,v)$-th entry, 
\[ e_{(a,v)}^T \left( (U^t x_0) \circ \comp{(U^t x_0)} \right),\]
gives the probability that the quantum walk lands on the particular arc $(a,v)$. We will call the following sum the \textsl{search probability} at time $t$:
\[ \sum_{v\sim a} e_{(a,v)}^T \left( (U^t x_0) \circ \comp{(U^t x_0)} \right).\]
An important question in quantum walk search is to determine, for certain family of graphs, the optimal time $t$ at which the search probability is sufficiently large.

In this paper, we study a related concept called the \textsl{average search probability}. While 
\[(U^t x_0) \circ \comp{(U^t x_0)}\]
does not converge, its time average converges, and the limit can be expressed using the spectral idempotents of $U$. This was observed in Aharonov, Ambainis, Kempe, and Vazirani \cite{Aharonov2000}, and we will state the result following the notation in Godsil and Zhan \cite{Godsil2017}.

\begin{lemma}
Let $U$ be a unitary transition matrix and let $x_0$ be a unit vector. The time-averaged probability distribution
\[\frac{1}{T}\sum_{t=0}^{T-1} (U^t x_0) \circ \comp{(U^t x_0)}\]
converges as $T$ goes to infinity. Moreover, if $F_r$ is the orthogonal projection onto the $r$-th eigenspace of $U$, then the limit is
\[\sum_r (F_r x_0) \circ \comp{(F_r x_0)}.\]
\end{lemma}

For the rest of the paper, unless otherwise specified, we will let $U$ and $x_0$ be
\[U= RCO_a, \quad \frac{1}{\sqrt{nk}}\one_n \otimes \one_k,\]
and let the spectral decomposition of $U$ be
\[U = \sum_r e^{i\theta_r} F_r.\]
 The \textsl{average search probability} of our quantum walk is
\[   \sum_r e_{(a,v)}^T \left( (F_r x_0) \circ \comp{(F_r x_0)} \right).\]


\section{Spectral decomposition}
In this section, we prove a spectral correspondence between $U$ and the weighted adjacency matrix of an augmented graph of $X$. 

We first note that, upon reordering the rows and columns,
\[CO_a = \pmat{
 I - \frac{2}{k}J & & &\\
 & \frac{2}{k}J-I & &\\
 & & \ddots &\\
 &&& \frac{2}{k}J-I}.\]
 Since this is a reflection, we can write it as twice a projection minus the identity. To be more precise, let $K$ be the $k\times (k-1)$ Vandermonde matrix:
\begin{equation}\label{eqn:K}
K = \pmat{
1&1 & 1 &\cdots &1\\
 e^{2\pi i/k}&  e^{4\pi i/k} & e^{6\pi i/k} & \cdots & e^{2(k-1)\pi i/k}\\
 \vdots&\vdots & \vdots &\vdots& \vdots\\
 e^{2(k-1)\pi i/k}& e^{4(k-1)\pi i/k}& e^{6(k-1)\pi i/k}&\cdots&e^{2(k-1)^2\pi i/k} }.
 \end{equation}
 Then $KK^*/k$ is the projection onto the orthogonal complement of $\mathrm{span}\{\one_k\}$. Thus, the first block of $CO_a$ can be rewritten as
 \[I - \frac{2}{k}J = \frac{2}{k} KK^*- I.\]
Now let $N$ be the block diagonal matrix
\begin{equation}\label{eqn:N}
N = \pmat{
K & & &&\\
 & \one_k & &&\\
 & & \one_k &&\\
 &&& \ddots&\\
 &&&&\one_k }
 \end{equation}
 with exactly $n$ blocks in the diagonal.  It follows from our discussion that 
 \[CO_a = \frac{2}{k} NN^* - I_{nk},\]
 that is, $CO_a$ is a reflection about the column space of $N$.
 
 On the other hand, the arc-reversal matrix $R$ is also a reflection. Thus, we may apply the following lemma, with
 \[P=\frac{1}{2}(R+I),\quad L=\frac{1}{\sqrt{k}}N,\]
  to determine the spectrum of $U$.
 
 \begin{lemma}\cite[Ch 2]{Zhan2018} \label{lem:spd}
 Let $P$ and $Q$ be two projections, and write $Q=LL^*$ for some matrix $L$ with orthonormal columns. Let
 \[U = (2P-I) (2Q-I).\]
Then the eigenspaces of $U$ are given as follows.
 \begin{enumerate}[(i)]
 \item The $1$-eigenspace of $U$ is the direct sum
 \[(\col(P)\cap \col(Q)) \oplus (\ker(P)\cap\ker(Q)).\]
 \item The $(-1)$-eigenspace of $U$ is the direct sum
 \[(\col(P)\cap \ker(Q)) \oplus (\ker(P)\cap\col(Q)).\]
 \item The remaining eigenspaces of $U$ are completely determined by the eigenspaces of $L^*(2P-I)L$. To be more specific, let $\lambda$ be an eigenvalue of $L^*(2P-I)L$ that lies strictly between $-1$ and $1$, and write $\lambda =\cos(\theta)$ for some $\theta\in\re$. The map
 \[z\mapsto ((\cos(\theta)+1)I - (e^{i\theta}+1)P)Lz\]
 is an isomorphism from the $\lambda$-eigenspace of $L^*(2P-I)L$ to the $e^{i\theta}$-eigenspace of $U$, and the map
 \[z\mapsto ((\cos(\theta)+1)I - (e^{-i\theta}+1)P)Lz\]
 is an isomorphism from the $\lambda$-eigenspace of $L^*(2P-I)L$ to the $e^{-i\theta}$-eigenspace of $U$.
 \end{enumerate}
 \end{lemma}
 
The above lemma shows that the eigenspaces of $U$ are largely determined by those of $L^*(2P-I)L$. For our quantum walk,
 \[P=\frac{1}{2}(R+I),\quad L=\frac{1}{\sqrt{k}}N,\]
 and so
 \[L^*(2P-I)L = \frac{1}{k}\pmat{ 0 & K^* & 0\\
 K &\multicolumn{2}{c}{
 \multirow{2}{*}{$A(X\backslash a)$}} \\
 0&},\]
 where the rows of $K$ are indexed by the neighbors of $a$. The right hand side, up to a scalar, is the Hermitian adjacency matrix of a weighted graph, obtained from $X$ by cloning the marked vertex $a$ and assigning appropriate $k$-th roots of unity to each arc. 


We can say more about the eigenspaces of $U$.

\begin{theorem}\label{thm:spd}
Let $X$ be a $k$-regular graph on $n$ vertices. Let $a$ be the marked vertex, and $R$, $C$ and $O_a$  the corresponding arc-reversal matrix, coin matrix and oracle. Let $K$ and $N$ be defined as in \eqref{eqn:K} and \eqref{eqn:N}. Then the eigenspaces of the transition matrix, $U=RCO_a$, satisfy the following. 
\begin{enumerate}[(i)]
\item There is a one-to-one correspondence between each conjugate pair of non-real eigenvalues $e^{\pm i\theta}$ of $U$, and each eigenvalue $\lambda$ of 
\begin{equation}\label{eqn:Ahat}
\widetilde{A}=\pmat{ 0 & K^* & 0\\
   K &\multicolumn{2}{c}{
   \multirow{2}{*}{$A(X\backslash a)$}} \\
   0&}
 \end{equation}
that lies strictly between $-k$ and $k$, and they are related by 
\[\lambda = k\cos(\theta).\]
Moreover, if $F_{\theta}$ is the projection onto the $e^{i\theta}$-eigenspace of $U$, and $\widehat{E_{\lambda}}$ is the projection onto the $\lambda$-eigenspace of $\widetilde{A}$, then
\[F_{\theta}=\frac{1}{2k\sin^2(\theta)} (N - e^{i\theta} RN) \widehat{E_{\lambda}}(N-e^{i\theta} RN)^*.\] 
\item If $X$ is $2$-connected, then the $1$-eigenspace of $U$ is orthogonal to $\mathrm{span}\{\one_{nk}\}$.
\end{enumerate}
\end{theorem}
\proof
We first prove (i). The eigenvalue correspondence follows from Lemma \ref{lem:spd}. To compute $F_{\theta}$, it suffices to show that when
 \[P=\frac{1}{2}(R+I),\quad L=\frac{1}{\sqrt{k}}N,\]
 the matrix
\[((\cos(\theta)+1)I - (e^{i\theta}+1)P)L\]
is a scalar multiple of 
\[N-e^{i\theta} RN,\]
and that for any $\lambda$-eigenvector $z$ of $\widetilde{A}$,
\[\norm{(N-e^{i\theta} RN)z}^2 = 2k\sin^2(\theta)\norm {z}^2.\]
Indeed,
\begin{align*}
((\cos(\theta)+1)I - (e^{i\theta}+1)P)L&=\frac{1}{2\sqrt{k}}\left((e^{-i\theta}+1) I -(e^{i\theta}+1)R\right)N\\
&=\frac{1}{2\sqrt{k} (e^{-i\theta}+1)} (I-e^{i\theta}R)N,
\end{align*}
and 
\begin{align*}
\norm{(N-e^{i\theta} RN)z}^2&=
z^*N^*(I-e^{-i\theta}R)(I-e^{i\theta}R)Nz\\
&=2z^*(N^*N -\cos(\theta) N^*RN)z\\
&=2 z^*(kI_{n+k-2} - \cos(\theta)\widetilde{A})z\\
&= 2(k-k\cos^2(\theta))\norm{z}^2\\
&=2k\sin^2(\theta). 
\end{align*}

We now prove (ii). By Lemma \ref{lem:spd}, the $1$-eigenspace of $U$ is
 \[(\col(R+I)\cap \col(NN^*)) \oplus (\ker(R+I)\cap\ker(NN^*)).\]
We will show that the second space in the direct sum is orthogonal to $\mathrm{span}\{\one_{nk}\}$, while the first space is trivial.

Since $R$ reverses each arc, the column space of $R+I$ consists of vectors that are constant on each pair of opposite arcs $(u,v)$ and $(v,u)$; in particular, $\one_{nk}$ lies in this space. Thus, $\ker(R+I)$ is orthogonal to $\mathrm{span}\{\one_{nk}\}$. 

On the other hand, any vector $x$ lying in the column space of $NN^*$ sums to zero over the outgoing arcs of $a$, 
and is constant on the outgoing arcs of any unmarked vertex; that is, 
\begin{equation}\label{eqn:zerosum}
\sum_{v\sim a} x_{(a,v)}=0,
\end{equation}
and for each $u\ne a$,
\begin{equation}\label{eqn:cons}
x_{(u,v)} = x_{(u,w)},  \quad \forall v\sim u, w\sim u.
\end{equation}
If, in addition, $x$ lives in the column space of $R+I$, then for each arc $(u,v)$,
\[x_{(u,v)} = x_{(v,u)}.\]
Now take any two neighbors $b$ and $c$ of $a$. Since $X$ is $2$-connected, there is a cycle containing edges $\{a, b\}$ and $\{a, c\}$, say
\[a, b, \seq{v}{1}{2}{\ell}, c, a.\]
Then
\[x_{(a,b)} = x_{(b,a)} = x_{(b, v_1)} = x_{(v_1, b)} = \cdots = x_{(c,a)} = x_{(a,c)}.\]
As our choice of $b$ and $c$ are arbitrary, this means that $x$ is constant on all outgoing arcs of $a$, which, together with \eqref{eqn:zerosum} and \eqref{eqn:cons}, implies that $x$ is the zero vector.
\qed

\section{Equitable partitions \label{sec:equitable}}

Given a graph $X$, a partition $\pi$ of its vertex set
\[\pi=\{\seq{C}{0}{1}{d}\}\]
is \textsl{equitable} if for any $i$ and $j$, there is a constant $c_{ij}$ such that every vertex in $C_i$ has $c_{ij}$ neighbors in $C_j$. Equivalently, if $P$ denotes the characteristic matrix of $\pi$, then $\pi$ is equitable if and only if there is some matrix $B$, called the \textsl{quotient matrix}, such that $A(X)P=PB$.

The following are standard results on equitable partitions; for more background, see Godsil \cite[Ch 5]{Godsil1993}.

\begin{lemma}\cite[Ch 5]{Godsil1993} \label{lem:equitable}
Let $\pi$ be an equitable partition of $X$. Let $P$ be the characteristic matrix, and $B$ the quotient matrix. 
\begin{enumerate}[(i)]
\item If $Bx = \lambda x$, then $A(X)Px = \lambda Px$.
\item If $A(X)y = \lambda y$, then $y^TPB= \lambda y^T P$.
\end{enumerate}
In particular, the eigenvectors of $A(X)$ either sum to zero over each cell of $\pi$, or are constant on each cell of $\pi$.
\end{lemma}

Clearly, every equitable partition $\pi$ for $X$ gives rise to an equitable partition $\pi\backslash \{C_i\}$ for the subgraph $X\backslash C_i$, and the quotient matrix can be obtained from $B$ by deleting the row and the column indexed by $C_i$. This leads to the following observation.

\begin{lemma}\label{lem:rowsum}
Let $\pi$ be an equitable partition of $X$, with $C_i$ as one of its classes. For any integer $m\ge 0$, the vector $A(X\backslash C_i)^m \one$ is constant on the cells of $\pi\backslash \{C_i\}$.
\end{lemma}
\proof
Let $P$ be the characteristic matrix of $\pi$. Since $\pi$ is equitable, there is some matrix $B$ such that
\[A(X) P = P B.\]
Now let $Q$ be the matrix obtained from $P$ by deleting the column indexed by $C_i$, and $F$ the matrix obtained from $B$ by deleting the row and the column indexed by $C_i$. Then 
\[A(X\backslash C_i) Q =Q F.\]
Thus for any positive integer $m$,
\[A(X\backslash C_i)^m Q =Q F^m,\]
and so
\[A(X\backslash C_i)^m\one =A(X\backslash a)^mQ\one = QF^m\one,\]
which lies in the column space of $Q$.
\qed

The next result is a consequence of Lemma \ref{lem:rowsum}.

\begin{corollary}\label{cor:powerA}
Let $X$ be a $k$-regular graph on $n$ vertices. Suppose $X$ has an equitable partition, where the singleton $\{a\}$ and the neighborhood $N(a)$ are two of its classes. Let $K$ and $\widetilde{A}$ be defined as in \eqref{eqn:K} and \eqref{eqn:Ahat}. For any integer $m\ge 0$,
\[\widetilde{A}^m\pmat{0_{k-1}\\\one_{n-1} } = \pmat{0_{k-1}\\ A(X\backslash a)^m \one}.\]
\end{corollary}
\proof
Clearly, 
\[\widetilde{A}^0\pmat{0_{k-1}\\\one_{n-1} } = \pmat{0_{k-1}\\ A(X\backslash a)^0 \one}.\]
 Suppose for some integer $m\ge 1$,
 \[\widetilde{A}^{m-1}\pmat{0_{k-1}\\\one_{n-1} } = \pmat{0_{k-1}\\ A(X\backslash a)^{m-1} \one}.\]
 Then
 \begin{align*}
 \widetilde{A}^m\pmat{0_{k-1}\\\one_{n-1} } &= \widetilde{A}\pmat{0_{k-1}\\ A(X\backslash a)^{m-1} \one}\\
 &=\pmat{ 0 & K^* & 0\\
    K &\multicolumn{2}{c}{
    \multirow{2}{*}{$A(X\backslash a)$}} \\
    0&}\pmat{0_{k-1}\\\one_{n-1} } \\
 &=\pmat{ \pmat{K^*& 0} A(X\backslash a)^{m-1}\one\\ A(X\backslash a)^m\one }.
 \end{align*}
By Lemma \ref{lem:rowsum}, $A(X\backslash a)^{m-1}\one$ is constant on $N(a)$, and since $K^*\one =0$,  the top block in the last vector vanishes. 
\qed

With $\widetilde{A}$ and $A(X\backslash a)$ defined above,  we now show a relation between their eigenspaces. To start, we cite a well-known result in linear algebra.

\begin{lemma}\label{lem:poly}
Let $M$ be a normal matrix, with distinct eigenvalues $\mu_r$ and eigenprojections $E_r$. Let
\[p_r(x) = \prod_{s\ne r}(x-\mu_s).\]
Then 
\[p_r(M) = p_r(\mu_r) E_r.\]
In particular, $E_r$ is a polynomial in $M$.
\end{lemma}

An eigenvalue of a graph is called a \textsl{main eigenvalue} if its eigenspace is not orthogonal to $\mathrm{span}\{\one\}$.

\begin{corollary}
Let $X$ be a $k$-regular graph on $n$ vertices. Suppose $X$ has an equitable partition, where the singleton $\{a\}$ and the neighborhood $N(a)$ are two of its classes. Let $K$ and $\widetilde{A}$ be defined as in \eqref{eqn:K} and \eqref{eqn:Ahat}. Then every main eigenvalue of $X\backslash a$ is an eigenvalue of $\widetilde{A}$. Moreover, if $\lambda$ is an eigenvalue of $\widetilde{A}$ with eigenprojection $\widetilde{E}_{\lambda}$, then 
\[ \widetilde{E}_{\lambda} \pmat{0_{k-1}\\\one_{n-1} } = \pmat{0_{k-1}\\ E_{\lambda} \one}\]
if $\lambda$ is also an eigenvalue of $X\backslash a$ with eigenprojection $E_{\lambda}$, and 
\[ \widetilde{E}_{\lambda} \pmat{0_{k-1}\\\one_{n-1} } = 0\]
otherwise.
\end{corollary}
\proof
Let $\mu$ be a main eigenvalue of $X$ with eigenprojection $E_{\mu}$. By Lemma \ref{lem:poly}, 
\[E_{\mu} = q (A(X\backslash a))\]
for some polynomial $q(x)$. It follows from Corollary \ref{cor:powerA} that 
\[q(\widetilde{A}) \pmat{0_{k-1}\\ \one_{n-1}}= \pmat{0_{k-1}\\ E_{\mu} \one }.\]
Now multiply both sides by $\widetilde{A}$. We have
\[\widetilde{A}q(\widetilde{A}) \pmat{0_{k-1}\\ \one_{n-1}}=\widetilde{A}\pmat{0_{k-1}\\ E_{\mu}\one }.\]
On the other hand,
\[\widetilde{A}q(\widetilde{A}) \pmat{0_{k-1}\\ \one_{n-1}}=\pmat{0_{k-1}\\ A(X\backslash a) q(A(X\backslash a)\one}=\pmat{0_{k-1}\\ A(X\backslash a) E_{\mu}\one}= \mu \pmat{0_{k-1} \\ E_{\mu}\one}.\]
As $E_{\mu}\one \ne 0$, 
\[\pmat{0_{k-1} \\ E_{\mu}\one}\]
is an eigenvector for $\mu{A}$ with eigenvalue $\mu$. This proves the first statement.

To see the second statement, let $\lambda$ be an eigenvalue of $\widetilde{A}$ with eigenprojection $\widetilde{E}_{\lambda}$. Again, Lemma \ref{lem:poly} tells us that $\widetilde{E}_{\lambda}$ is a polynomial in $\widetilde{A}$; moreover, this polynomial sends $\lambda$ to $1$, and all other eigenvalues of $\widetilde{A}$ to $0$. Now, since each main eigenvalues of $X\backslash a$ is an eigenvalue of $\widetilde{A}$, this polynomial also sends $A(X\backslash a)$ to its $\lambda$-eigenspace if $\lambda$ is indeed an eigenvalue of $A(X\backslash a)$, and to $0$ otherwise.
\qed


\section{Distance regular graphs \label{sec:drg}}
One special type of equitable partitions arise in distance regular graphs. A graph is called \textsl{distance regular} if for any two vertices $u$ and $v$ at distance $m$, the number of vertices at distance $i$ from $u$ and distance $j$ from $v$ is a constant which depends only on $i$, $j$ and $m$. As a result, for any vertex $a$, the distance partition relative to $a$ is equitable, and the quotient matrix is tridiagonal:
\begin{equation}\label{eqn:B}
B=\pmat{
a_0 & b_0 && & &\\
c_1 & a_1 & b_1& &&\\
& \ddots & \ddots & \ddots &&\\
&& c_{d-1} &a_{d-1} & b_{d-1}\\
&& &c_{d} & a_{d}}.
\end{equation}
Moreover, upon permuting the rows and columns, this quotient matrix does not depend on the choice of $a$. The list of parameters
\[\{\seq{b}{0}{1}{d-1}; \seq{c}{1}{2}{d}\}\]
is called the \textsl{intersection array} of the distance regular graph. 

We cite some basic results on the intersection array. For references, see the book by Brouwer, Cohen and Neumaier \cite[Sec 4.1]{Brouwer1989}, or the survey by van Dam, Koolen and Tanaka \cite{Dam2018}.

\begin{lemma}\label{lem:drg}
Let $X$ be a distance regular graph, with valency $k$ and intersection array $\{\seq{b}{0}{1}{d-1}; \seq{c}{1}{2}{d}\}$. The following hold.
\begin{enumerate}[(i)]
\item $a_0=0$, $b_0=k$ and $c_1=1$.
\item For each $i$,
\[a_i = k - b_i - c_i,\]
with the convention that $c_0=b_d=0$.
\item For any vertex $a$, the number $k_i$ of vertices at distance $i$ from $a$ satisfies
\[b_ik_i = c_{i+1} k_{i+1}.\]
\end{enumerate}
\end{lemma}

We now prove a lemma on a submatrix of the Laplacian matrix of a distance regular graph; this turns out useful when we study the average search probability in Section \ref{sec:avgprobdrg}.

\begin{lemma}\label{lem:inreasing}
Let $X$ be a distance regular graph on $n$ vertices with valency $k$. Let $L(X)$ be the Laplacian matrix of $X$. For any vertex $a$, let $L(x)\backslash a$ denote the matrix obtained from $L(X)$ by removing the $a$-th row and the $a$-th column. Then the entries in $(L(X)\backslash a)^{-1} \one$ increase in distance from $a$. Moreover, for any neighbor $v$ of $a$ in $X$,
\[e_v^T (L(X)\backslash a)^{-1} \one=\frac{n-1}{k}.\]
\end{lemma}
\proof
It is a well-known result that $L(X)\backslash a$ is invertible. Let 
\[y =(L(X)\backslash a )^{-1} \one.\]
As $X$ is distance regular, the distance partition relative to $a$ gives an equitable partition for $X\backslash a$, and so by Lemma \ref{lem:rowsum} and the fact that $(L(X)\backslash a )^{-1}$ is a polynomial in $L(X)\backslash a$, the vector $y$ is constant on the cells of this partition. That is, there are integers $\seq{z}{1}{2}{d}$ such that, for each vertex $u$ at distance $i$ from $a$ in $X$, the entry $y_u$ equals $z_i$.

Now we solve for $y$ in 
\begin{equation}\label{eqn:solve}
(L(X)\backslash a )y = \one.
\end{equation}
Pre-multiply both sides by $\one^T$, and we get
\[\pmat{\one_k\\ 0}^T y = \one^T (L(X)\backslash a) y = \one^T \one =n-1.\]
Therefore,
\[z_1 = \frac{n-1}{k}.\]
Moreover, expanding \eqref{eqn:solve} gives
\[kz_1 - a_1z_1 - b_1z_2=1,\]
and, assuming $z_{d+1}=0$,
\[kz_i - c_i z_{i-1} - a_i z_i -b_i z_{i+1}=1,\quad i=2,\cdots d.\]
Since $k=a_i+b_i+c_i$, this is equivalent to
\[b_1(z_2-z_1)= z_1 - 1\]
and
\[b_i(z_{i+1}-z_i) +1= c_i(z_i-z_{i-1}),\quad i=2\cdots, d\]
Let $k_i$ denote the number of vertices at distance $i$ from $a$ in $X$. Using Lemma \ref{lem:drg} and induction, we see that
\[z_{i+1}-z_{i} = \frac{k_{i+1}+\cdots+k_d}{k_i b_i}>0.\]
Hence $y$ increases in distance from $a$.
\qed

Finally, we cite a result on vertex-connectivity, due to Brouwer and Koolen \cite{Brouwer2009}.

\begin{theorem}\cite{Brouwer2009}.
The vertex-connectivity of a distance regular graph equals its valency.
\end{theorem}

\section{Average search probabilities for distance regular graphs \label{sec:avgprobdrg}}

We now consider the search problem on distance regular graphs. Let $X$ be a distance regular graph on $n$ vertices with valency $k\ge 2$. Let $a$ be the marked vertex. Let $U=RCO_a$ be as defined in Section \ref{sec:model}, and $F_{\theta}$ the projection onto its $e^{i\theta}$-eigenspace of $U$. Recall that the average search probability is 
\[\frac{1}{nk}\sum_r e_{(a,v)}^T \left( (F_{\theta} \one ) \circ \comp{(F_{\theta} \one)} \right).\]
We say an eigenvalue $e^{i\theta}$ \textsl{contributes to} search if $F_{\theta}\one \ne 0$.

\begin{lemma}\label{lem:eav}
The eigenvalues of $U$ that contribute to search are $-1$ and $e^{\pm i\theta}$, where $\theta = \arccos(\lambda/k)$ for some eigenvalue $\lambda$ of $X\backslash a$. Moreover, if $E_{\lambda}$ is the projection onto the $\lambda$-eigenspace of $X\backslash a$, then for any neighbor $v$ of $a$ in $X$,
\[e_{(a,v)}^T F_{\theta}\one  = \frac{1}{2\sin^2(\theta)} (1-e^{i\theta}) e_v^T E_{\lambda} \one.\]
\end{lemma}
\proof
We apply Theorem \ref{thm:spd}. As $X$ is $2$-connected, the eigenvalue $1$ of $U$ does not contribute to search. Let $e^{i\theta}$ be a non-real eigenvalue of $U$. Let $\widetilde{A}$ be defined as in \eqref{eqn:Ahat}. Then $\lambda=k\cos(\theta)$ is an eigenvalue of $\widetilde{A}$, and 
\[F_{\theta}=\frac{1}{2k\sin^2(\theta)} (N - e^{i\theta} RN) \widetilde{E_{\lambda}}(N-e^{i\theta} RN)^*,\]
where $\widetilde{E}_{\lambda}$ is the projection onto the $\lambda$-eigenspace of $\widetilde{A}$, and $N$ is defined as in \eqref{eqn:N}. To see whether $e^{i\theta}$ contributes to search, we multiply both sides by $\one$. Since
\[N^*R\one = N^*\one=k\pmat{0_{k-1}\\ \one_{n-1}},\]
we get
\[F_{\theta}\one = \frac{1}{2\sin^2(\theta)} (1-e^{-i\theta})(N-e^{i\theta} RN) \widetilde{E_{\lambda}}\pmat{0_{k-1}\\ \one_{n-1}},\]
which, by Corollary \ref{cor:powerA}, further reduces to
\[\frac{1}{2\sin^2(\theta)} (1-e^{-i\theta})(N-e^{i\theta} RN) \pmat{0_{k-1}\\ E_{\lambda}\one},\]
where $E_{\lambda}$ is the projection to the $\lambda$-eigenspace of $X\backslash a$. Now, suppose $v$ is the $j$-th neighbor of $a$ in $X$. We have
\[e_{(a,v)}^T (N-e^{i\theta} RN) \pmat{0_{k-1}\\ E_{\lambda}\one_{n-1}}= (e_j^T \pmat{K & 0}  - e^{i\theta} e_v^T)\pmat{0_{k-1}\\ E_{\lambda}\one} =-e^{i\theta} e_v^T E_{\lambda} \one,\]
from which the statement follows.
\qed

The above lemma allows us to express the average search probability completely in terms of the spectrum of $X\backslash a$. 

\begin{theorem}\label{thm:avgsearch}
Let the spectral decomposition of $X\backslash a$ be
\[A(X\backslash a) = \sum_{\lambda} \lambda E_{\lambda}.\]
Let $v$ be any neighbor of $a$ in $X$. Then the average search probability of the quantum walk on $X$ is given by
\[\frac{1}{n} \sum_{\lambda} \frac{k^3}{(k-\lambda)(k+\lambda)^2} (e_v^T E_{\lambda}\one)^2+\frac{1}{n} \left(1-\sum_{\lambda} \frac{k}{k+\lambda} e_v^T E_{\lambda} \one \right)^2.\]
\end{theorem}
\proof
The first term follows from Lemma \ref{lem:eav} by taking the square of $e_{(a,v)}^T F_{\theta}\one$. For the second term, notice that 
\[F_{\pi} = I - \sum_{\theta \ne \pi } F_{\theta},\]
and so the $(-1)$-eigenprojection of $U$ is determined by the remaining eigenspaces of $U$. Moreover, since the non-real eigenvalues of $U$ come in conjugate pairs,
\[F_{\theta} = \comp{F_{-\theta}}.\]
This together with Lemma \ref{lem:eav} yields the second term in the statement.
\qed

\section{Tridiagonal matrices and orthogonal polynomials \label{sec:tridiag}} 
In this section, we discuss the connection between eigenvectors of a generic tridiagonal matrix, and eigenvectors of a graph whose quotient matrix coincides with this tridiagonal matrix. This will eventually provide us tools to study the limit of the average search probability on a family of distance regular graphs.

Let $T$ be an $m\times m$ tridiagonal matrix,
\begin{equation}\label{eqn:T}
T=\pmat{
\alpha_0 & \beta_0 && & &\\
\gamma_1 & \alpha_1 & \beta_1& &&\\
& \ddots & \ddots & \ddots &&\\
&& \gamma_{m-1} &\alpha_{m-1} & \beta_{m-1}\\
&& &\gamma_m & \alpha_m},
\end{equation}
where the diagonal entries $\alpha_i$ are non-negative, and the off-diagonal entries $\beta_i$, $\gamma_i$ are positive. Let $T_i$ denote the leading $i\times i$ principal matrix of $T$. Then $T$ defines a sequence of monic polynomials,
\begin{align*}
p_0(x)&=1\\
p_1(x)&= \det(xI-T_1)\\
p_2(x)&=\det(xI-T_2)\\
&\vdots\\
p_m(x)&=\det(xI-T_m)\\
p_{m+1}(x)&=\det(xI-T_{m+1}),
\end{align*}
called the \textsl{orthogonal polynomials} associated with $T$. Note that $p_{m+1}(x)$ is the characteristic polynomial of $T$.

We cite some standard results on orthogonal polynomials. The proofs can be found in Szego \cite{Szego1939} or Chihara \cite{Chihara2011}.

\begin{theorem}
For each $i$, the roots of $p_i(x)$ are real and simple. Moreover, the roots of $p_i(x)$ interlace those of $p_{i+1}(x)$.
\end{theorem}

The orthogonal polynomials associated with $T$ determine its eigenvectors. 



\begin{theorem}\label{thm:zlambda}
Let $\lambda$ be a root of $p_m(x)$. Let $z_{\lambda}$ be the vector defined by
\[z_{\lambda} = \pmat{p_0(\lambda)\\ p_1(\lambda)/\beta_0\\ \vdots \\p_m(\lambda)/(\beta_0\cdots \beta_{m-1})}.\]
Then $z_{\lambda}$ is an eigenvector for $T$ with eigenvalue $\lambda$.
\end{theorem}

In Section \ref{sec:drg}, we saw that each distance regular graph has a quotient matrix of the form \eqref{eqn:B}, and the entries are determined by the intersection array 
\[\{\seq{b}{0}{1}{d-1}; \seq{c}{1}{2}{d}\}.\]
In fact, we may consider the ``reversed" tridiagonal matrix as well:
\begin{equation}\label{eqn:S}
S=\pmat{
a_d & c_d && & &\\
b_{d-1} & a_{d-1} & c_{d-1}& &&\\
& \ddots & \ddots & \ddots &&\\
&& b_1 & a_1 & c_1\\
&& &b_0 & a_0},
\end{equation}
and, as we will see later, this form becomes handy when we relate the eigenvectors of a distance regular graph to that of its vertex-deleted subgraphs. The orthogonal polynomial associated with $S$ are known as the \textsl{dual orthogonal polynomials} to those associated with \eqref{eqn:B}, and has been studied in, for example,  Vinet and Zhedanov \cite{Vinet2004}.

From now on, assume $X$ is a distance regular graph of diameter $d\ge 2$ on $n$ vertices with valency $k$. Fix a vertex $a$, and let $k_i$ be the number of vertices in $X$ at distance $i$ from $a$. Let 
\[P = \pmat{
\one_{k_d} & & &&\\
 & \one_{k_{d-1}} & &&\\
 & & \ddots &&\\
 &&& \ddots \one_{k_1}&\\
 &&&&\one_{k_0} }.\]
 Then $P$ is the characteristic matrix of the ``reversed" distance partition, whose quotient matrix is defined in \eqref{eqn:S}. As before, let $S_i$ denote the leading $i\times i$ principal matrix of $S$, and let $\{q_i(x)\}$ be the orthogonal polynomials associated with $S$, that is,
 \[q_i(x) = \det(xI-S_i).\]
We also define a new sequence of polynomials by $u_0(x)=1$ and 
 \[u_i(x) = \frac{q_i(x)}{c_d c_{d-1} \cdots c_{d-i+1}}, \quad i=2,3, \cdots d.\]

\begin{lemma}\label{lem:bound}
Let $\lambda$ be the largest eigenvalue of $X\backslash a$ and $E_{\lambda}$ the corresponding eigenprojection. For any neighbor $v$ of $a$ in $X$, we have
\[e_v^T E_{\lambda} \one\ge \frac{n-1}{n} (u_{d-1}(\lambda))^2= \frac{n-1}{n} \left(\frac{q_{d-1}(\lambda)}{q_{d-1}(k)}\right)^2.\]
\end{lemma}
\proof
The largest eigenvalue of a graph is a main eigenvalue, and so by Lemma \ref{lem:equitable}, it must be an eigenvalue of its quotient matrix. Consider $X$ first. Since it is $k$-regular and $S$ is its quotient matrix, the vector
\[z_k=\pmat{u_0(k)\\ u_1(k)\\ \vdots \\u_{d-1}(k)\\u_d(k)}\]
coincides with the all-ones vector, and $Pz_k$ is an eigenvector for $X$ with eigenvalue $k$. We have
\[n=\norm{Pz_k}^2 = z_k^T P^TP z_k.\]
On the other hand, the submatrix $S_d$ is the quotient matrix of $X\backslash a$. Thus,
\[q_d(\lambda)=0,\]
and
\[ y_{\lambda}=\pmat{u_0(\lambda)\\ u_1(\lambda)\\ \vdots \\u_{d-1}(\lambda)}\]
is an eigenvector for $S_d$ with eigenvalue $\lambda$. Let $Q$ be the matrix obtained from $P$ by deleting the last column. Then $Qy_{\lambda}$ is an eigenvector for $X\backslash a$ with eigenvalue $\lambda$. Moreover, by interlacing, for each $i$ we have
\[u_i(k)\ge u_i(\lambda)>0,\]
and so
\begin{align*}
\norm{Q y_{\lambda}}^2 &= y_{\lambda}^T Q^T Qy_{\lambda}\\
&=y_{\lambda}^T Q^T Qy_{\lambda} + u_d(\lambda)\\
&\le z_k^T P^TP z_k\\
&=n.
\end{align*}
Therefore,
\begin{align*}
e_v^T E_{\lambda} \one &= \frac{1}{\norm{Qy_{\lambda}}^2} \ip{e_{d-1}}{Q y_{\lambda}} \ip{\one}{Q y_{\lambda}}\\
&\ge \frac{1}{n} u_{d-1}(\lambda)( k_d u_0(\lambda) + k_{d-1} u_1(\lambda) + \cdots + k_1u_{d-1}(\lambda)).
\end{align*}
Now, since $y_{\lambda}$ is an eigenvector for $S_d$, we have a recurrence relation:
\[b_{d-i} u_{i-1}(\lambda) + a_{d-i} u_i(\lambda)+ c_{d-i}u_{i+1}(\lambda) = \lambda u_i(\lambda),\]
from which and $\lambda<k$ we get
\[b_{d-i}(u_{i-1}(\lambda)-u_i(\lambda))< c_{d-i} (u_i(\lambda) - u_{i+1}(\lambda))\]
for $i\ge 1$, and 
\[u_1(\lambda)<u_0(\lambda).\]
Thus, $u_i(\lambda)\ge u_{d-1}(\lambda)$ for all $i\le d-1$, and so
\[e_v^T E_{\lambda} \one \ge \frac{1}{n} (u_{d-1}(\lambda))^2 (k_d+\cdots + k_1) = \frac{n-1}{n} (u_{d-1}(\lambda))^2.\]
Finally, observe that $u_{d-1}(k)=1$.
\qed

Now we show that the largest eigenvalue $\lambda$ of $X\backslash a$ gets closer to $k$ as $k$ increases. A nice result of Eldridge, Belkin and Wang \cite{Eldridge2018}, derived from Courant-Fischer-Weyl min-max principle, turns out to be useful.

\begin{theorem}\cite{Eldridge2018}\label{thm:perturb}
Let $A$ and $B$ be two $n\times n$ real symmetric matrix. Let $\mu_1\ge \mu_2 \ge \cdots \ge \mu_n$ be the eigenvalues of $A$, with corresponding eigenvectors $\seq{x}{1}{2}{n}$. Let $\lambda_1\ge \lambda_2 \cdots \ge \lambda_n$ be the eigenvalues of $A+B$. Given $j\in\{1,2,\cdots,n\}$, suppose 
\[h\ge \abs{\ip{x}{Bx}}\]
for all unit vectors $x$ in $\mathrm{span}\{\seq{x}{1}{2}{j}\}$. Then for $i=1,2,\cdots,j$,
\[\lambda_i \ge \mu_j -h.\]
\end{theorem}

Applying this to the symmetrized quotient matrix

\begin{equation}\label{eqn:Shat}
\widehat{S}=\pmat{
a_d & \sqrt{b_{d-1}c_d} && & &\\
\sqrt{b_{d-1}c_d} & a_{d-1} & \sqrt{b_{d-2}c_{d-1}}& &&\\
& \ddots & \ddots & \ddots &&\\
&& \sqrt{b_1c_2} & a_1 & \sqrt{b_0c_1}\\
&& &\sqrt{b_0c_1} & a_0},
\end{equation}
we find a lower bound for the largest eigenvalue of $X\backslash a$.

\begin{lemma}\label{lem:labound}
Let $\lambda$ be the largest eigenvalue of $X\backslash a$. Then
\[\lambda \ge k -\frac{2k}{n}.\]
\end{lemma}
\proof
It is not hard to see that $\widehat{S}$ is similar to $S$, and 
\[A P (P^TP)^{-1/2} =P (P^TP)^{-1/2} \widehat{S}. \]
Therefore, $k$ is the largest eigenvalue of $\widehat{S}$, with unit eigenvector 
\[x=\frac{1}{\sqrt{n}} \pmat{\sqrt{k_d}\\ \sqrt{k_{d-1}} \\ \vdots \\ \sqrt{k} \\ 1}.\]
We consider the symmetrized quotient matrix of $X\backslash a$, which is the leading $d\times d$ principal submatrix of $\widehat{S}$. Write 
\[\pmat{
a_d & \sqrt{b_{d-1}c_d} && & &\\
\sqrt{b_{d-1}c_d} & a_{d-1} & \sqrt{b_{d-2}c_{d-1}}& &&\\
& \ddots & \ddots & \ddots &&\\
&& \sqrt{b_1c_2} & a_1 & 0\\
&& &0 & 0} = \widehat{S}-\pmat{
0 & 0 && & &\\
0 & 0& 0& &&\\
& \ddots & \ddots & \ddots &&\\
&& 0& 0 & \sqrt{k}\\
&& &\sqrt{k} & 0}.\]
Clearly, $\lambda$ is the largest eigenvalue of the left hand side. Since 
\[\pmat{\sqrt{k} & 1} \pmat{0 & \sqrt{k} \\ \sqrt{k} & 0} \pmat{\sqrt{k}\\ 1}=2k,\]
by Theorem \ref{thm:perturb}, $\lambda\ge k -2k/n$.
\qed

\section{The limit of the average search probability}
We would like to know how the average search probability, on a family of distance regular graphs, behaves as the valency grows. As before, let $X$ be a distance regular graph on $n$ vertices, with valency $k\ge 2$ and a fixed diameter $d$. Let $a$ be the marked vertex, and $v$ any neighbor of $a$. Let the spectral decomposition of $X\backslash a$ be
\[A(X\backslash a) = \sum_{\lambda} \lambda E_{\lambda}.\]
 Theorem \ref{thm:avgsearch} says that the average search probability on $X$ splits into two terms:
\begin{align*}
s_1 &= \frac{1}{n} \sum_{\lambda} \frac{k^3}{(k-\lambda)(k+\lambda)^2} (e_v^T E_{\lambda}\one)^2\\
s_2 &= \frac{1}{n}\left(1-\sum_{\lambda} \frac{k}{k+\lambda} e_v^T E_{\lambda} \one \right)^2
\end{align*}
We show that $s_2$ vanishes as the graph gets larger.

\begin{lemma}\label{lem:s2}
As $n$ goes to infinity, $s_2$ approaches $0$.
\end{lemma}
\proof
The eigenvalues $\lambda$ of $X\backslash a$ interlace those of $X$. In particular, $-k<\lambda <k$. Hence
\[\frac{k}{k+\lambda}>\frac{1}{2},\]
and so
\[e_v^T \left(\sum_{\lambda} \frac{k}{k+\lambda} E_{\lambda}\right) \one >\frac{1}{2} e_v^T \left(\sum_{\lambda} E_{\lambda} \right)\one = \frac{1}{2}.\]
The result now follows from the squeeze theorem.
\qed

Our next observation gives a lower bound for $s_1$.

\begin{lemma}\label{lem:s1}
Let $\alpha$ be the smallest entry in $(L(X)\backslash a)^{-1}\one$. Then
\[s_1 > \frac{1}{4}\frac{n-1}{n} \sum_{\lambda} (E_{\lambda} J E_{\lambda})_{vv}.\]
\end{lemma}
\proof
We can rewrite $s_1$ as
\[s_1 = \frac{1}{n} e_v^T \left(\sum_r \frac{k^3} {(k-\lambda)(k+\lambda)^2} E_{\lambda} J E_{\lambda} \right) e_v,\]
which, by interlacing, is greater than
\[\frac{1}{4n} e_v^T \left(\sum_{\lambda} \frac{k}{k-\lambda} E_{\lambda} J E_{\lambda} \right)e_v.\]
Since
\[L(X)\backslash a = kI - A(X\backslash a) = \sum_{\lambda} (k-\lambda) E_{\lambda},\]
we have
\begin{align*}
\sum_{\lambda} \frac{k}{k-\lambda} E_{\lambda} J E_{\lambda} &=k \left(\sum_{\lambda} \frac{1}{k-\lambda} E_{\lambda}\right) \left( \sum_{\lambda} E_{\lambda} J E_{\lambda}\right)\\
&=k (L(X)\backslash a)^{-1} \left( \sum_{\lambda} E_{\lambda} J E_{\lambda}\right)\\
&=k \sum_{\lambda} E_{\lambda}(L(X)\backslash a)^{-1} J E_{\lambda}.
\end{align*}
The inequality now follows  Lemma \ref{lem:inreasing}.
\qed

For readers who are familiar with continuous-time quantum walks, we notice that 
\[\frac{1}{n-1} \sum_{\lambda} E_{\lambda} J E_{\lambda}\]
is the \textsl{average state} of the continuous-time quantum walk on $X\backslash a$, with $\one/\sqrt{n-1}$ as the initial state. This reveals an interesting connection between the discrete-time quantum walk on a graph and the continuous-time quantum walk on its vertex-deleted subgraph. More discussion on average states can be found in Coutinho, Godsil, Guo, and Zhan\cite{Coutinho2018}.

We now find the limit for the search probabilities on complete graphs.

\begin{theorem}
The average search probability on $K_n$ approaches $1/4$ as $n$ goes to infinity.
\end{theorem}
\proof
The vertex deleted subgraph $K_{n-1}$ has two spectral idempotents:
\[E_{n-2} = \frac{1}{n-1} J,\quad E_{-1}=-\frac{1}{n-1}J.\]
As they are orthogonal to each other, we have
\[\sum_{\lambda} E_{\lambda} J E_{\lambda} = E_{n-2} J E_{n-2} = J.\]
Hence the average probability converges to $1$.
\qed

Any family of \textsl{strongly regular graphs}, that is, distance regular graphs of diameter two, also enjoy this probability.

\begin{theorem}
The average search probability on a strongly regular approaches $1/4$ as the valency goes to infinity.
\end{theorem}
\proof
Let $X$ be a strongly regular graph with intersection array $\{k, b_1; 1, c_2\}$. Clearly, $b_1=k-a_1-1$, and $a_2=k-c_2$. The main eigenvalues of $X\backslash a$ are eigenvalues of 
\[S=\pmat{a_1 & k-a_1-1\\ c_2 & k-c_2}.\]
Solving $Sz = \lambda z$ yields 
\[\lambda = \frac{1}{2} (k+a-c) \pm \frac{1}{2} \sqrt{(k-a+c)^2-ac}\]
and 
\[z = \frac{1}{\sqrt{(k-a-1)c + (\lambda-a)^2}} \pmat{\sqrt{(k-a-1)c}/\sqrt{k}\\ (\lambda-a)/\sqrt{a-k-1}}.\]
Thus, if $E_{\lambda}$ is the $\lambda$-eigenprojection for $X\backslash a$ and $v$ is any neighbor of $a$ in $X$, then by our discussion above, the average search probability is at least a quarter of 
\[ \sum_{\lambda} (E_{\lambda} J E_{\lambda})_{vv}=\frac{(k-a_1+c_2)^2-4c_2}{(k-a_1+c_2)^2 - 4c_2 - 2(k-a_1-1)}.\]
Since $a_1, c_2\le k$, the right hand side approaches $1$ as $k$ goes to infinity.
\qed

For distance regular graphs with larger diameters, we need the results from Section \ref{sec:tridiag} to bound the limit.

\begin{theorem}\label{thm:criterion}
Let $d\ge 3$ be a fixed positive integer. Let $X$ be a distance regular graph of diameter $d$. Assume the number of vertices $n$, the valency $k$, and the intersection numbers $a_i, b_i, c_i$ are all functions in a parameter  $\tau$, and $k(\tau)$ increases in $\tau$. Further assume 
\[\lim_{\tau=0}\frac{k(\tau)^{d-1}}{c_2(\tau) \cdots c_d(\tau) n(\tau)}=0.\]
Let $a$ be any vertex of $X$. As the the valency goes to infinity, the average search probability on $X$ approaches $1/4$.
\end{theorem}
\proof
Let the quotient matrix $S$ of $X$ be as defined in \eqref{eqn:S}. Let $\{q_i(x)\}$ be the associated sequence of orthogonal polynomials.  Let $n$ be the number of vertices in $X$, and $v$ any neighbor of $a$. Let $\lambda$ be the largest eigenvalue of $X\backslash a$, with eigenprojection $E_{\lambda}$. Lemma \ref{lem:bound} tells us that
\[e_v^T E_{\lambda} \one \ge \frac{n-1}{n} \left(\frac{q_{d-1}(\lambda)}{q_{d-1}(k)}\right)^2.\]
Since $a_i, b_i, c_i \le k$, by an inductive argument, the characteristic polynomial of any $\ell \times \ell$ principal submatrix of $S$,
\[x^{\ell} + \alpha_1 x^{\ell-1} +\alpha_2 x^{\ell-2}+ \cdots + \alpha_{\ell},\]
satisfies $\alpha_i = O(k^i)$. Thus
\begin{align*}
q_{d-1}\left( k -\frac{2k}{n}\right) &= q_{d-1}(k) + q'_{d-1}(k)\left(-\frac{2k}{n}\right)+q''_{d-1}(k)\left(-\frac{2k}{n}\right)^2+\cdots\\
&=q_{d-1}(k) + \frac{O(k^{d-1})}{n} +  \frac{O(k^{d-1})}{n^2}+\cdots + \frac{O(k^{d-1})}{n^{d-1}}.
\end{align*}
The result now follows from $q_{d-1}(k)=c_2c_3\cdots c_d$.
\qed

A special case, with stronger assumptions, is given below.

\begin{theorem}\label{thm:criterion2}
Let $d\ge 3$ be a fixed positive integer. Let $X$ be a distance regular graph of diameter $d$. Assume the number of vertices $n$, the valency $k$, and the intersection numbers $a_i, b_i, c_i$ are all functions in a parameter $\tau$, and $k(\tau)$ increases in $\tau$. Further assume 
\[\lim_{\tau\to0}\frac{k(\tau)^{d-1}}{n(\tau)}=0.\]
Let $a$ be any vertex of $X$. As the the valency goes to infinity, the average search probability on $X$ approaches $1/4$.
\end{theorem}

We note that this criterion is met by many common families of distance regular graphs, including the Hamming graphs $H(d,\tau)$, the Johnson graphs $J(\tau, d)$, the Grassmann graphs $J_q(\tau, D)$ with a fixed $q$, and the dual polar graphs with a fixed $e$.

\section{Future work}
The average probability we studied lies in the following vector,
\[\sum_r  (F_r x_0) \circ \comp{(F_r x_0)} ,\]
which, as we have seen, is the limit of 
\[\frac{1}{T}\sum_{t=0}^{T-1} (U^t x_0) \circ \comp{(U^t x_0)}\]
as $T$ goes to infinity. For graphs with high average search probabilities, it will be helpful to determine the \textsl{mixing time} $M_{\epsilon}$, that is, the smallest $K$ such that for all $T>K$,
\[\norm{\frac{1}{T}\sum_{t=0}^{T-1} (U^t x_0) \circ \comp{(U^t x_0)}-\sum_r  (F_r x_0) \circ \comp{(F_r x_0)}}\le\epsilon.\]
This could potentially indicate that quantum search is fast on some graphs.


\section*{Acknowledgment}
This project is supported by the York Science Fellows program. The author would like to thank Ada Chan, Chris Godsil and Thomas Wong for helpful discussions.

\bibliographystyle{amsplain}
\bibliography{dqw.bib}

\end{document}